# Elimination of spiral waves in 2D excitable media
## by  feedback control

Yelena Smagina and Moshe Sheintuch,

*Department of Chemical Engineering,*

Technion, Haifa, Israel, 32000

## Abstract

Control strategy for suppression of spiral-wave in a  2-D model of an excitable  media  is developed with application to the cardiac system.  The controller  which  incorporates a finite number of actuators (electrodes) assures  the establishment of   traveling  plane behavior using  the following strategy in the context of cardiac system: a small current that is proportional to a discrepancy between  actual  and  assigned voltages at some sensor position is imposed  via the  actuators. We present a systematic methodology for controller design and for stability analysis  based on  an approximate  1-D model of  the wave  front. The validity  of control  is checked by numerical simulations. Control is successful  provided  a  sufficient  number  of  actuators  is  incorporated.  The proposed control method is similar in spirit to that the previously developed by the authors  for a 1-D model of  cardiac tissue.

## 1.  Introduction

Spiral  waves  are  observed  in  many   different  excitable  media  of   physical  and biological  systems [1,2]. In some cases, spiral patterns    are undesirable   and it is necessary to suppress them and reset the system to its usual    behavior. For example, normal heart activity operation requires  traveling plane waves which move with constant shape.  A sufficiently strong perturbation can cause traveling plane waves breakup and give rise to spiral waves which may be dangerous: if there exist defective region in a heart (due to myocardial infarction), spirals may be  trapped  by  defects,  and  cause permanent  tachycardia [3].   Usually  these  spiral  waves  drift  and   dissipate  at  tissue borders and  as a result  the heart returns to its normal rhythms.

The problem of spiral-wave control has attracted significant attention and several strategies were suggested [4-17]. Non-feedback control by inserting an external forcing signal in order to move the spiral wave in the desired direction was applied in several mainly numerical works [4-8]. Feedback control of various types were published in [9-17]. Several authors use global (averaged over some domain) feedback via electric field modulation [9-14]. Other researchers apply external electrochemical devices that compensate wave disruption by control of the difference between actual and set points [15-17]. For instance, in [15] a multiple point-sensor/actuator time-delay control via series of electrodes implanted on the heart was used to numerically test the control. In theoretical work [16] the feedback controller is constructed to drive the system having spiral waves back to the plane wave behavior. In numerical work [17] the feedback signal is recorded at a certain time and is returned into the system when it enters the recovery state. The main drawback of approaches [16], [17] is the need to apply measures and actuators over the whole domain.

In the suggested work we also develop a feedback control that suppresses a spiral wave and resets the traveling plane behavior. The controller consists from finite number of electrodes implanted on the heart and realizes the following strategy: small current that is proportional to the discrepancy between actual and assigned voltages at some sensor position, affects the system via finite number of actuators and forcing system to return to the plane waves. This type of control is robust: it effective, in spite of the perturbations in the system parameters. For example, we check the effectiveness of the controller in terms of heterogeneity when one parameter of the model varies in some domain. The proposed control method is similar in spirit to the previously developed by the authors in [18]. The considered approach is more practical because it uses a smaller number of actuators and sensors: in previous elaborated control methods [15,16]) control is applied by sensors and actuators over the whole domain.

The main challenges of this control problem stem from the moving nature of patterns, while the actuators are likely to be stationary, from the 2-D nature of the problem and from the control law necessary. While we use here sensors that measure local properties (transmembrane potential, $V$) and actuators that reply locally, other approaches that use global sensing (e.g., average $V$ over the tissue) or actuating may be theoretically more efficient, but we are not aware of such devices for measuring cardiac properties. The control should keep the front planar: we show that this can be achieved by exciting the system when it lags behind the required front.

The structure of this paper is as follows; In the next section we outline the model, explain the excitable nature of the system with the activator and inhibitor as its variables and state the control problem. Section 3 describes the structure of the elaborated control law. In Section 4 we conduct an approximate analysis using an approximate model that keeps planar front motion in a bistable system (by setting the inhibitor value to be fixed); the analysis proves the control efficiency under approximation and reveals characteristic length scales for sensing and actuating. Finally the effectiveness of the control law is numerically tested in sections 5 with various densities of actuators: and in section 6 with heterogeneity in the model.

## 2. Statement

In this paper we consider a 2D popular version of the classical Fitzhugh-Nagumo two-variable model that contains most of the qualitative dynamic features of biological phenomena relevant to the formation of spiral and irregular waves. The aim of this work is to design a feedback control law which keeps the system at its traveling plane wave solution and corrects deviations if necessary.

We consider the general 2D model described in [16]

$$\frac{\partial V}{\partial t} = D(\frac{\partial^2 V}{\partial x^2} + \frac{\partial^2 V}{\partial y^2}) + f(V) - W + I(x, y, t) \tag{1}$$

$$\frac{\partial W}{\partial t} = \varepsilon(\beta V - \gamma W + \delta) \equiv \varepsilon g(V, W) \tag{2}$$

$$f(V) = V(\alpha - V)(V - 1) \tag{3}$$

where $V$ is the fast variable, an activator that is related to transmembrane potential in cardiac tissue, and $W$ is the slow variable, an inhibitor that expresses the effect of ion channels on membrane current in cardiac tissue; $x \in [0, L_x]$, $y \in [0, L_y]$, $I(x, y, t)$ represents injected perturbation (control). We set $D = 1$ and the multiplier of $f(V)$ to be unity so that the length scales are made dimensionless. Other parameters $\varepsilon = 0.01, \alpha = 0.1, \quad \beta = 0.5, \quad \gamma = 1$ and $\delta = 0$ are chosen to assure excitable media behavior [16]: (see Fig.1a).

Here we use two types of boundary conditions imposed at the edges of the rectangular two-dimensional domain: either no-flux

$$\frac{\partial V}{\partial x}(0, y) = \frac{\partial V}{\partial x}(L_x, y) = 0, \quad \frac{\partial V}{\partial y}(x, 0) = \frac{\partial V}{\partial y}(x, L_y) = 0 \qquad (4a)$$

or periodic conditions in the y-direction

$$V(x, 0) = V(x, L_y), \quad \frac{\partial V}{\partial y}(x, 0) = \frac{\partial V}{\partial y}(x, L_y) \qquad (4b)$$

while maintaining no-flux boundary conditions in x-direction. This transforms the rectangular domain into a torus. Periodic boundary conditions make it possible to choose a much smaller length of the domain in the direction of moving plane waves than that required for the no-flux conditions.

The 1-D solution of (1-3) takes the form of a pulse that moves in this excitable media in constant form and velocity (Figure 1b). Note that the pulse combines two fronts (jump in $V$) that move at constant $W$ which adjusts very slowly ($\varepsilon << 1$). The ascending and decreasing fronts of the pulse move with equal speed. In a plane this solution takes the form of plane wave which may loose its stability. A circulating plane wave is produced by perturbing part of the domain to obtain a propagating plane wave: i.e. we simulate Eqn. (1)-(3) with no-flux boundary conditions and initial condition: $V|_{t=0} = 1$ for $0.245L_y \le y \le 0.26L_y, 0 \le x \le L_x$ and $V|_{t=0} = V_s$ for otherwise ; $W|_{t=0} = W_s$ where $V_s, W_s$ is the steady state of equations (Eqn.1,2): $V_s = 0$, $W_s = 0$. We obtain two plane waves that move in opposite directions. When one moving wave runs out of the domain, the remaining single moving wave is kept running by joining numerically two ends of the rectangular domain into a torus.

The steady plane solution $V_o = V_o(x, y, t)$ of Eqn.1-3 possesses a constant shape (see Fig.2a) and constant traveling velocity ($c_{oy}$) [16]. This traveling wave front is perpendicular to the y-axis. Spirals can be produced up by breaking up a part of the front (zeroing it): $V|_{pert} = 0$, $0.05L_x \le x \le 0.5L_x$, $0 \le y \le L_y$. The broken wave wanders around this forcing part and gives rise to formation of a spiral (see Fig.2b). This class of behavior is the source of cardiac dysfunctions such as arrhythmia.

Our aim is to design a feedback control which should be able to drive the system dynamics back to the plane wave behavior characterized by a nonzero velocity $c_y = c_{oy}$ at all points of plane wave front in y-direction and zero velocity in x-direction $c_x \approx 0$. When the plane wave is broken then the velocity in x-direction grows slowly but the equality $c_x << c_y$ is satisfied initially at every point $x \in [0, L_x]$. We will try establish plane wave behavior while initially $c_x << c_y$. From this point on we may apply velocity control in y-direction.

Below we keep in mind that $c_y = c_y(x,t)$ may vary along the wave (in x-direction).

Thus we formulate problem as follows

**Problem:** Find control current $I = I(x, y, t)$ such that $c_y(x,t) \to c_{oy}$ at all points of the plane wave front, i.e. control that stabilizes velocity at all points of the front plane wave $V$ solution against small perturbations and despite parameter uncertainty. That implies also determining the locations and numbers of sensors and actuators.

We note that although we focus on a FN model, the same control could be applied for a wide range of RD systems of relevance to physiological or chemical process. This control should be understood as external electro-mechanical device that are able to compensate the signal disruption.

### 3. Structure of control

We will find a feedback controller (regulator) having the structure

$$I(x, y, t) = \bar{c}\, \frac{\partial V}{\partial y}$$

by assigning variable $\bar{c}$ a value proportional to the difference $V - V_o$ at some assigned sensor position when the front travels through it (or $\bar{c}$ is function of $V - V_o$ at the several assigned sensor positions $i = 1, 2, \ldots \eta$ situated at front line $y = Z$ ) i.e.

$$\bar{c} \sim (V(x_i, Z, t) - V_o(x_i, Z, t))$$

with $i \in [1, \ldots, \eta]$. For small variation between $V$ and the desired value $V_o$ we may approximate

$$V(x_i, Z, t) - V_o(x_i, Z, t) \sim \left. \frac{\partial V}{\partial y} \right|_{y=Z} ((Z(x_i, t) - Z_o(x_i, t))$$

where $Z(x_i, t)$ and $Z_o(x_i, t)$ are front line positions of actual and desired moving planar solutions $V$ and $V_o$ at a some assigned location $x_i$.

To study stability of the system with control we consider the system in a coordinate moving in $y$ – direction : $(y, t) \to (\varsigma, t')$ where $\varsigma = y - c_{oy} t$, $t' = t$, $c_{oy}$ is a front velocity of unperturbated steady solution in $y$ -direction . After transformation

$$\frac{\partial V}{\partial y} = \frac{\partial V}{\partial \varsigma} \ , \quad \frac{\partial V}{\partial t} = \frac{\partial V}{\partial t'} - c_{oy} \frac{\partial V}{\partial \varsigma}$$

Eqn.1 becomes ( with $D=1$)

$$\frac{\partial V}{\partial t'} = (\frac{\partial^2 V}{\partial x^2} + c_{oy} \frac{\partial V}{\partial \varsigma} + \frac{\partial^2 V}{\partial \varsigma^2}) + f(V) - W + I(x, \varsigma, t') \qquad (1a)$$

Thus Eqns.(1a, 2,3) describe dynamics of plane wave $V$ in a moving coordinate. The control term $I(x, \varsigma, t)$ [1] is calculated here as $I(x, \varsigma, t) = \bar{c} \frac{\partial V}{\partial \varsigma}$ where the feedback $\bar{c} \sim Z(x, t) - Z_o$ where $Z_o$ is assigned constant front line in moving coordinate . At a some positions $x_i$ we have

$$\bar{c}(x_i, t) = -k(Z(x_i, t) - Z_o) \ , \ \ i \in [1, \ldots, \eta] \qquad (6)$$

In Eqn.6 $k > 0$ is the gain coefficient . Note that control (6) affects all points of the domain. In the next section we consider a control that affects only finite number of points of the domain.

Note that when a discrete version of the model is solved (as we do below) $Z_o$ is fixed and should be placed at a line with sensors. Also in a laboratory-frame the front moves through sensors and is sensed for short time intervals only.

### 4.Stability analysis using approximate model

---

[1] Further for simplicity we will use $t$

As explained earlier the pulse is composed of two fronts moving at the same velocity. Stability of planar fronts to transversal perturbations, and possible control remedies, were studied in [19] in the context of chemical reactors. To substantiate the structure of the control law let's consider the following approximations: When $V$ is fast while $W$ is the slow variables ( i.e. $\varepsilon << 1$ in Eqns. 1-3 ) we can obtain a first approximation of the velocity of a $V$ front for frozen $W$ solution. Let $Z(x,t)$ be the planar front position in moving coordinate[2] and $c_\infty = c_\infty(W_o(Z(x,t)),I)$ is the approximate velocity of a planar front in y-direction that depends on the local steady-state inhibitor ( $W_o$ ) and control variable ( $I$ ). The velocity is also curvature dependent and for Eqn.1 with diffusity scaled into the variables, the velocity of a low- curvature front is $c_y \cong c_\infty - Z_{xx}$ [19] [3]. Thus considering the position of an ascending ( $dV/d\varsigma > 0$ )[4] planar front $Z = Z(x,t)$ and assuming a small perturbation from the stationary planar front we approximate the velocity $c_y$ as [19]

$$-\frac{dZ}{dt} = c_y = c_\infty(W(Z(x,t),I)) - \frac{\partial^2 Z}{\partial x^2}$$

where the sign of $dZ/dt$ implies that the low state of the front is on the left and the high on the right and $c_\infty > 0$ implies expansion of latter.

Let $Z_o(x) = Z_o$, a constant value, be the assigned front line in a moving coordinate that is perpendicular to y-direction. The front propagation equation in moving coordinate linearized around $Z_o$ ($I=0$) becomes

$$-\frac{d\overline{Z}}{dt} \cong (\frac{\partial c_\infty}{\partial W})_f (\frac{\partial W}{\partial \varsigma})_f \overline{Z} - \frac{\partial^2 \overline{Z}}{\partial x^2} \tag{7}$$

where $\overline{Z} = Z(x,t) - Z_o$, and the subscript 'f' denotes that derivatives are estimated at front line. Since $\varepsilon(\beta V_o - \gamma W_o + \delta) = 0$ (Eqn.2) at steady state then we obtain the

---

[2] We propose that the ascending and descending fronts of the pulse move with equal speed.

[3] Recall that in moving coordinate $c_\infty = 0$ for planar front.

[4] Let's note that the theoretical results ( elaborated for the negative direction of the wave) are true for the wave moving in positive direction.

following approximate expression $(\partial W / \partial \varsigma)_f = \frac{\beta}{\gamma} (\partial V / \partial \varsigma)_f$ [5] for small deviations of

$V_o$ and $W_o$. As a result equation (7) becomes

$$\frac{d\overline{Z}}{dt} \cong -\frac{\beta}{\gamma} (\frac{\partial c_\infty}{\partial W})_f (\frac{\partial V}{\partial \varsigma})_f \overline{Z} + \frac{\partial^2 \overline{Z}}{\partial x^2} \qquad (8)$$

When $\partial c_\infty / \partial W < 0$ then the uncontrolled ($I=0$) system (8) may be unstable to transversal perturbations. We consider fronts that admits non-flux boundary conditions . Thus we will find $\overline{Z}$ as $\overline{Z} \sim a_n(t)\cos(n\pi x / L_x) + \vartheta$, $n = 1,2,...$ where $\vartheta$ is a certain constant value.

Substituting this $\overline{Z}$ in above equation we approximate the eigenvalues of Eqn.8 as

follows

$$\rho_n = -\frac{\beta}{\gamma} (\frac{\partial c_\infty}{\partial W})_f (\frac{\partial V}{\partial \varsigma})_f - (\frac{n\pi}{L_x})^2, \quad n = 1,2,... \qquad (9)$$

The most unstable eigenvalue of the open-loop equation Eqn. 8 is $\rho_1 = \chi - (\pi / L_x)^2$

where the constant $\chi = -\frac{\beta}{\gamma} (\frac{\partial c_\infty}{\partial W})_f (\frac{\partial V}{\partial \varsigma})_f$ is evaluated at the frozen profile (see Appendix).

For our parameters $\chi = 3.54 > 0$.

Now we demonstrate that adding a simple feedback control $I = I(x,t)$ of the form $I(x,t) = -\widetilde{k}\overline{Z}(x,t)$ in right hand side of Eqn.8, where $\overline{Z}(x,t)$ is deviation of front line $Z(x,t)$ from the set one $Z_o$ and $\widetilde{k}$ is some gain coefficient, is sufficient to stabilize the solution of the closed-loop system obtained. Hence, the planar front be maintained in the transversal directions at the set line front $Z_o$ . We will seek gain $\widetilde{k}$ as $\widetilde{k} = k(\frac{\partial V}{\partial \varsigma})_f$ thus

control takes the form

$$I(x,t) = -k(\frac{\partial V}{\partial \varsigma})_f \overline{Z}(x,t) \qquad (10)$$

Calculating the eigenvalues of the linearized closed-loop system (Eqn. 8 with control

(10)) :

---

[5] Here $\delta = 0$.

$$\frac{d\overline{Z}}{dt} = -\frac{\beta}{\gamma}(\frac{\partial c_\infty}{\partial W})_f (\frac{\partial V}{\partial \varsigma})_f \overline{Z} + \frac{\partial^2 \overline{Z}}{\partial x^2} - k(\frac{\partial V}{\partial \varsigma})_f \overline{Z}$$

gives the following eigenvalues

$$\rho_n = -\frac{\beta}{\gamma}(\frac{\partial c_\infty}{\partial W})_f (\frac{\partial V}{\partial \varsigma})_f - (\frac{n\pi}{L_x})^2 - k(\frac{\partial V}{\partial \varsigma})_f, \quad n = 1,2,... \tag{11}$$

Since $(\partial V/\partial \varsigma)_f > 0$ then unstable eigenvalue of the open-loop system may be shifted to the negative domain by imposing a control term $k(\frac{\partial V}{\partial \varsigma})_f$ with a sufficiently large gain $k$.

Now consider the realistic point-sensor control that should maintain the plane front in the transversal directions at the set front line $Z_o$ (in approximate model in moving coordinate)

$$I(x,t) = -k(\frac{\partial V}{\partial \varsigma})_f \sum_{d=1}^{\mu} \overline{Z}(x_e,t)\, \psi_d(x) \tag{12}$$

where $x_e, \ e \in [1,...,\eta_x], \quad \eta_x < \mu^6$ are the sensor positions along the front line, $\psi_d(x) = \delta(x - x_d^*)$ are actuator functions that realize the pointwise control and are applied at points $x_d^*, \ d = 1,2,...,\mu$. Here $\delta$ is 1D Dirac delta function.

Control (12) has a general structure and is very difficult to realize in the considered application (heart) because it responds to many discrepancies $\overline{Z}(x_e,t) = Z(x_e,t) - Z_o$. To simplify its realization we suggest to measure only a single discrepancy $\overline{Z}(x_l,t)$, $l \in [1,\eta_x]$. Thus control (12) becomes the special a form

$$I(x,t) = -k(\frac{\partial V}{\partial \varsigma})_f ((Z(x_l,t) - Z_o)\sum_{d=1}^{\mu} \delta(x - x_d^*) \tag{13}$$

Now we will study the influence of the actuators number ($\mu$) to assure linear stability of the closed-loop system. Let us expand $\overline{Z}(x,t)$ and $\overline{Z}(x_l,t)$ as

$$\overline{Z}(x,t) = \sum_n a_n(t)\phi_n(x), \quad \overline{Z}(x_l,t) = \sum_n a_n(t)\phi_n(x_l) \tag{14}$$

where $\phi_n(x)$, the eigenfunctions of the problem,

---

[6]We propose that $\eta_x < \mu$ : number sensors less than number of actuators.

$$\frac{\partial^2 \phi(x)}{\partial x^2} = -\rho\phi(x), \qquad \varphi_x(0) = \varphi_x(L_x) = 0 \qquad (15)$$

are

$$\phi_n(x) = \begin{cases} 1/\sqrt{L_x}, & n = 1 \\ \dfrac{\sqrt{2}}{\sqrt{L_x}}\cos\dfrac{\pi(n-1)x}{L_x}, & n = 2,\ldots \end{cases}$$

(16)

In (15) eigenvalues $\rho$ are

$$\rho_n = n^2\pi^2/L_x^2, \qquad n = 1,2,\ldots \qquad (17)$$

We lump linearized equation (8) with control (13) by Galerkin method to obtain the following realization

$$\dot{a}_n = [-\frac{\beta}{\gamma}(\frac{\partial c_\infty}{\partial W})_f(\frac{\partial V}{\partial y})_f - \rho_n]a_n - k(\frac{\partial V}{\partial y})_f b_n \sum_j a_j \phi_j(x_d), \ n = 1,2,\ldots \qquad (18)$$

where $\rho_n$ are the eigenvalues of the above linear operator and

$$b_1 = \sum_{d=1}^{\mu}(1/\sqrt{L_x}), \qquad b_n = (\sqrt{2}/\sqrt{L_x})\sum_{d=1}^{\mu}\cos\frac{n\pi}{L_x}x_d, \ n = 2,\ldots \qquad (19)$$

In vector-matrix notation Eq.(18) becomes a linear infinite-dimensional dynamical system with scalar input $v$ and output $w$

$$\dot{a} = Aa + Bv, \qquad\qquad w = Ha \qquad (20)$$

closed by the finite-dimensional linear output feedback

$$v = -k(\frac{\partial V}{\partial \varsigma})_f w \qquad (21)$$

In (20) $a=a(t)$ is the infinite dimensional vector; $v$ and $w$ are scalars; matrix $A$ has structure

$$A = diag(\chi - \frac{\pi^2}{L_x^2}, \chi - \frac{4\pi^2}{L_x^2}, \ldots, \chi - \frac{n^2\pi^2}{L_x^2}, \ldots) \qquad (22)$$

where $\chi = -\frac{\gamma}{\beta}(\frac{\partial c_\infty}{\partial W})_f(\frac{\partial V}{\partial y})_f$. Elements of infinite-dimensional column vector $B$ calculated in (19) and infinite-dimensional row vector $H$ has elements

$$H_j = \phi_j(x_l), \ j = 1,2,\ldots, \ l \in [1,\eta] \qquad (23)$$

Below we will study a truncated (finite-dimensional) approximation of Eqn.(20) with truncated order *N*. An approximate estimate of truncation order *N* may be obtained by conventional methods.

Thus, the problem may be restated as follows: For the linearized truncated ODEs (20) it is necessary to find output feedback control (21) that stabilizes the closed-loop system. In other words we need to choose sensor positions $x_e$, $e \in [1, \ldots, \eta]$ and actuator positions $x_d^*$, $d = 1, 2, \ldots, \mu$ to form vectors *B* and *H* that assure solvability of the problem.

Solvability conditions here coincide with solvability conditions of widely known problem of output stabilization of a multivariable system with single input and output. They are reduced to the following conditions on vectors B and H (for detail see series works [19]-[21] with references to works with the concept of system zeros).

1. $\det(HB) \neq 0$,

2. roots of the equation $B(sI - A)^{-1} H = 0$ are negative,

where matrix A satisfies (22) and elements of *B* and *H* are calculated according formulas (19), (23) respectively.

These conditions help to fulfil rough estimation of sensor and actuator positions and to evaluate number of actuators.

**Remark 1.** Control (13) is continuous in time. The realistic discrete-time analogy of (13) is $I(x, t_j) = -k(\frac{\partial V}{\partial \varsigma})_f ((Z(x_l, t_j) - Z_o) \sum_{d=1}^{\mu} \delta(x - x_d^*)$. We may show that negativeness of eigenvalues of matrix $A - kBH$ guarantees effectiveness of discrete-time control. Indeed, in this case Eqn. (20) becomes a hybrid one $\dot{a} = Aa + Bv(t_j)$, $w(t_j) = Ha(t_j)$. To simplify the asymptotic stability analysis of the following closed-loop hybrid ODEs

$$\dot{a}(t) = Aa(t) - kBHa(t_j) \tag{24}$$

we sample the continuous-time part of above equation to obtain the following discrete–time closed-loop equation

$$a(t_{j+1}) = G(t_{j+1}, t_j) a(t_j), \quad j = 1, 2, \ldots$$

where the matrix $G(t_{j+1}, t_j)$, calculated as $G = \exp((A - kBH)(t_{j+1} - t_j))$, is the transition matrix. The above closed-loop discrete-time system is asymptotically stable if eigenvalues of the correspondent dynamic matrix $G(t_{j+1}, t_j) = \exp((A - kBH)(t_{j+1} - t_j))$ are within unit circle. This is satisfied if eigenvalues of the dynamic matrix of continuous time system $A - kBH$ are negative.

In the original 2D model if control (13) is implemented via actuators located on the lines $y_q^*$, $q = 1, \ldots, \mu_y$ (in y-direction) at the positions $x_d^*$, $d = 1, \ldots, \mu$ (in x-direction), takes the following structure

$$I(x, y, t) = -k(\frac{\partial V}{\partial y})_f ((Z(x_l, t) - Z_o(t)) \sum_{d=1}^{\mu} \sum_{q}^{\mu_y} \delta(x - x_d^*, y - y_q^*), \quad l \in [1, \ldots, \eta_x] \qquad (25)$$

where the subscript '$f$' denotes that derivatives are estimated at front line Z and $\delta(x - x_d^*, y - y_q^*)$ is 2D Dirac delta function : $\delta(x, y) = 0$ for $x \neq 0$, $y \neq 0$ and $\iint \delta(x, y) dx dy = 1$ for $x = x_d^*, y = y_q^*$.

**General scheme of control realization:** Let plane wave $V(x, y, t)$ moves in y-direction and its front line is be continuously measured at $\eta_x$ sensors situated at positions $x_i = x_1 + \Delta(i-1)$, $i = 1, \ldots, \eta_x$ in x-direction and at positions $y_1, \ldots, y_{\eta_y}$ in y-direction with $d_j = y_j - y_{j-1}$. When a single point $x_l$, $l \in [1, \ldots, \eta_x]$ of the front line reaches the assigned line $\tilde{y}_e$, $e \in [1, \ldots, \eta_y]$ then this moment is fixed as $t_j$ and the difference $\bar{c}(t_j) = d_j - d_{set}$ is used for control updating where $d_{set} = c_{set}(t_j - t_{j-1})$ with $c_{set}$ is assigned velocity. In notation (25) values $d_j$, $d_{set}$ correspond $Z$ and $Z_o$ respectively.

Updated control which operated via actuators situated at positions $(x_d^*)$, $d = 1, \ldots, \mu$ in x-direction and at lines $y_q^*$, $q \in [1, \ldots, \mu_y]$ in y-direction is given by

$$I(x, y, t_j) = -k(\frac{\partial V}{\partial y})_{\tilde{y}_l} ((Z(x_l, t_j) - Z_o) \sum_{d=1}^{\mu} \sum_{q=1}^{\mu_y} \delta(x - x_d^*, y - y_q^*), \quad l \in [1, \ldots, \eta_x], \qquad (26)$$

## 5. Simulations.

Now we demonstrate the effectiveness of control law (25).   All numerical simulations are performed in a two-dimensional $N_x \times N_y$ lattice with $N_x = 200$ and $N_y = 400$. We use a time step dt =0.5 and a space step *dx=dy=1.0*.

The effect of idealized controller (25) that uses large number of points of the space for realization is presented in Figure 3. Since spiral wave velocity is curvature dependent and is less than the velocity of plane wave, then a smaller number of sensors ($\eta_x = 2$) situated at distance $L_x / 3$, $2L_x / 3$ is enough for control realization.

By simulation we discover that we may reduce the number of actuators to a grid of 25 in x-direction ($\mu_x = 25$ at points 1,5,..,100) by 6 ($\mu_y = 6$ at points 110,140,…,260) in y-direction (see Fig. 4 for spatiotemporal patterns).

Increasing the number of actuators improves the response time of the system. So, doubling the number points in x-direction ($\mu_x = 50$) gives a significant decrease the time to returning to travel waves (see Figure 5).

## 6. Control in a heterogeneous medium

The ionic heterogeneity is considered to be one of the main factors underlying the initiation of spiral waves [22]. The aim of this section is to demonstrate that even when ionic heterogeneity initiates spiral waves in the simplest FitzHugh-Nagumo model (1)-(3), the suggested control can reset these waves back to the planar traveling waves form.

We add heterogeneity to model (1)-(3) by setting parameter $\alpha = 0.21$ in the domain $0.01L_x \le x \le 0.5L_x$, $0.275L_y \le y \le 0.775L_y$ ( $i \in [2, 100]$, $j \in [110, 310]$), while keeping the α=0.1 in the rest of the domain.. Also we propose that this heterogeneity applies during restricted time interval (that is typical for cardiac tissue parameters that change during infarction [23]). The above heterogeneity introduced at t=200 induces spirals solution (Fig.7) while setting control (25) under heterogeneity (Figure 8) shows the dynamics returning the planar wave behavior.

## Conclusion

The problems of pattern formation and selection are central to nonlinear science. In this paper we use control to suppress spiral wave patterns and encourage plane patterns in a 2-D model. Transitions from plane wave patterns to spiral ones may lead to serious cardiac diseases. Therefore, the problem of spiral wave suppression in order to reset the excitable media to a travelling plane patterns is of extreme importance in cardiology.

Now, in both nonlinear science and cardiac physiology fields there is a growing efforts to develop low-amplitude control methods. Our theoretical research is based on simplest 2-D model, of an excitable media. This model incorporates a pair of coupled nonlinear parabolic PDEs with polynomial source functions (reaction-diffusion PDEs). The parameters chosen assure that the system is locally stable but excitable. The low-amplitude control strategy is based on injecting small currents in response to discrepancies between actual and assigned voltages at some sensor positions.

We obtain stability conditions of the closed-loop system in terms of the linearized approximate 1-D model that describes dynamics of the moving planar front.

In contrast to similar studies (see [15,16]), a major advantage of the proposed method is the use of a relatively small number of actuators and sensors. In the cited works sensors and actuators are applied over the whole domain. The drawback of the proposed control is using an additional parameter $(\partial V / \partial y)_f$. Since the configuration of the pulse profile is kept at all times then the value $(\partial V / \partial y)_f$ may be evaluated theoretically. In any case, estimation of the actual front line may not be so easy to implement in practical situations. We hope to develop modified techniques that utilize the advantage of our approach and overcome its disadvantages.

## Acknowledgements

Work supported by the US-Israel Binational Science Foundation. MS is a member of the Minerva Center of Nonlinear Dynamics of Complex Physics. YM is supported by the Center for Absorption in Science, Ministry of Immigrant Absorption, State of Israel. We are grateful Dr. Olga Nekhamkina for useful discussion.

**APPENDIX**

**Calculating the right hand term in Eqn.8**

We calculate the term $\frac{\beta}{\gamma}(\frac{\partial c_\infty}{\partial W})_f(\frac{\partial V}{\partial y})_f$ by approximation of $c_\infty(W)$. Planar front velocity $(c_\infty)$ of the system $\frac{\partial V}{\partial t} = (\frac{\partial^2 V}{\partial x^2} + \frac{\partial^2 V}{\partial y^2}) + f(V) - W$ with set $W = W_o = \frac{\beta}{\gamma}V_o$

at the steady state value is calculated as [24]

$$c_\infty = \frac{1}{\sqrt{2}}(V_+ + V_- - 2V_i) \qquad (A1)$$

where $V_+$, $V_-$, $V_i$ are roots of nonlinear equation $V[(\alpha - V)(V - 1)] - W_o = 0$. To evaluate these roots we use linear Taylor's approximation.

To this point we denote $f(V) = V[(\alpha - V)(V - 1)] = -V^3 + V^2(\alpha + 1) - \alpha V$ and

$f_1(V) = V[(\alpha - V)(V - 1)] - W_o = -V^3 + V^2(\alpha + 1) - \alpha V - W_o$. Considering function

$f_1(V) = f(V) - W_o$ with $|W_o| << 1$ as a perturbation function we may write for its root

$V^* + \lambda$ where $V^*$ is corresponding root of $f(V) = 0$ and $|\lambda| << 1$ is the following

linear Taylor's approximation

$$f_1(V^* + \lambda) \cong f_1(V^*) + \frac{\partial f_1}{\partial V}\Big|_{V^*} \cdot \lambda = f(V^*) - W_o + (-3V^{*2} + 2V^*(\alpha + 1) - \alpha)\lambda = 0$$

Since $f(V^*) = 0$ then we obtain the following approximation for $\lambda$

$$\lambda = W_o / (-3V^{*2} + 2V^*(\alpha + 1) - \alpha) \qquad (A2)$$

Substituting $V_+ = 1$, $V_- = 0$, $V_i = 0.1$ instead of $V^*$ in (A2) we find $\lambda_+ = -W_o / 0.9$,

$\lambda_- = -W_o / 0.1$, $\lambda_i = W_o / 0.09$. Then, roots obtained $V_+ + \lambda_+$, $V_- + \lambda_-$ and $V_i + \lambda_i$ may

be used for calculating $c_\infty = -\frac{33.3W_0}{\sqrt{2}}$ and $(\frac{\partial c_\infty}{\partial W})_f = -33.3/\sqrt{2}$.

The derivative $(\frac{\partial V}{\partial y})_f$ may be evaluate by differential of analytical steady-state

solution of PDE (1a), (2) which takes the form $V_o = \sqrt{1 - \beta}\tanh(y - Z_f)\sqrt{0.5(1 - \beta)}$

[25]. We obtain $(\frac{\partial V}{\partial y})_f = (1 - \beta)/\sqrt{2}$ and $\frac{\beta}{\gamma}(\frac{\partial V}{\partial y})_f = \beta(1 - \beta)/(\gamma\sqrt{2}) = 0.177$.

More accurate result may be obtained numerically from the slope $(\frac{\partial W_o}{\partial \varsigma})_f \cong \frac{\beta}{\gamma}(\frac{\partial V}{\partial y})_f$ : e.g. in Fig. 1b we have $(\partial W_o / \partial \varsigma)_f = 0.15$ . Thus

$$\frac{\beta}{\gamma}(\frac{\partial c_\infty}{\partial W})_f (\frac{\partial V}{\partial y})_f = (\frac{\partial c_\infty}{\partial W_o})_f (\frac{\partial W_o}{\partial \varsigma})_f \cong -3.54 .$$

## References

[1] M.C. Cross and P.C. Hohenberg, Rev. Mod. Phys. **65**, 851(1993).

[2] D. J. Christini and L. Glass, Chaos, **12**, 732 (2002).

[3] A.T. Winfee, Science, **266**, 1003 (1994).

[4] V.N. Biktashev, and A.V. Holden, Chaos, **8**, 48 (1998).

[5] R.M. Mantel, and D. Barkley, Phys. Rev. E, **54**, 4791 (1996) .

[6] H. Zang , B. Hu, and G. Hu, Phys. Rev. E, **68**, 026134 (2003).

[7] H. Pueba , R. Martin, J. Alvarez-Ramirez , and R. Aguilar-Lopez , Chaos, Solutions and Fractals, **39,** 971 (2009).

[8] P-Y. Wang , P. Xie, and H-W Yin, Chin. Phys. Soc. **12,** 674 (2003).

[9] S.Grill , V.S. Zykov, and S.C. Muller, Phys. Rev. Lett, **75**, 3368 (1995).

[10] V.S. Zykov , A.S. Mikhailov, and S.C. Muller, Phys. Rev. Lett, **78**, 3398 (1998).

[11] G.Yuan , Xu Aiguo, G. Wang , and S. Chen, Europhysics Letters, **90**, 1 (2010).

[12] O –U. Kheowan, V.S. Zykov , and S.C. Muller, Physical Chemistry Chemical, **4,** 1334 (2002).

[13] O-U. Kheowan, and S.C. Muller, Applied Mathematics and Computation, **164**, 373 (2005).

[14] V.S. Zykov , and H. Engel, Physica D, **199**, 243 (2004).

[15] W.J. Rappel, F. Fenton, and A.Karma, Phys. Rev. Lett, **83**, 456 (1999).

[16] C. Vilas, M.R. Garcia, J.R. Banga, and A.A. Alonso, Physica D, **237**, 2353 (2008).

[17] G.Y.Yuan, S.G. Chen, and S.P. Yang, Eur. Phys. J. **B, 58,** 331 (2007).

[18] Y. Smagina, and M. Sheintuch. Proc. 18[th] World Congress IFAC, Milan, Italy, August 28-September 2, 5379 (2011).

[19] M. Sheintuch, Ye. Smagina, and O. Nekhamkina, Phys. Rev. E **66** , 066213

(2002).

[20] Ye Smagina,  M.Sheintuch, , Preprints of 5[th] IFAC Symposium on Robust
     Control Design,  Toulouse, France, July 5-7, 2006.

[21]  Ye. Smagina, O. Nekhamkina. and  M. Sheintuch, Journal of Process Control,
      **16**, (2006).

[22]  A. Defauw, P. Dawyndt and A.V. Panfilov, Phys. Rev. E, **88,** 062703-1 (2013).

[23]  R.R. Aliev and A.V.Panfilov,  J. Theor. Biol **181,** 33  (1996).

[24] A.S.  Mikhailov, *Foundation of Synergetics 1. Distributed Active Systems*
     (Springer-Verlag, Berlin, 1994).

[25]  M.Sheintuch, Ye. Smagina,  and  O. Nekhamkina,   Ind.  Eng. Chem. Res. **41,**
      2136 (2002) .

### Figures

**Figure1.** The nullcurves  *f(V) –W=0*  and  *g(V,W)=0*  of  the uncontrolled  system
(Eqn.1-3) ( $\varepsilon = 0.01, \alpha = 0.1,\ \ \beta = 0.5,\ \gamma = 1\ ,\delta = 0\ $ and $V_s = 0,\ W_s = 0$ ) **(a)**   and  the
corresponding 1-D solution  **(b)** .

**Figure 2.**  Spatiotemporal gray-scale planar **(a)**  or spiral waves **(b)** in the ($y,x$) plane of
the  uncontrolled  dimensionless  transmembrane     potential   $V(x, y,t)$   (Eqns.1-3,
$N_x = 200$ points,  $N_y = 400$ points,  parameters as in Fig.1); plane waves  break  at t=7  by
setting the domain   $0.05L_x \le x \le 0.5L_x,\ 0 \le y \le L_y$ to $V = 0$. Figures shows  snapshots at
times  t=100, t=200,….

**Figure 3.** Testing  effectiveness of idealized controller (25) to  suppress spiral waves  and
maintain planar fronts. The controller uses 7  fixed equally spaced   sensors in x-direction
( $y = 0.735L_y$ )  and   a large number of actuators ( a grid of 100 actuators in x-direction
and 200 actuators in y-direction ; k=0.001). Control start  at t=10. Here $c_{set} = 0.5$, other
parameters are as in Figure 1.  Figure shows  snapshots at times  t=200, t=400,….

**Figure 4.** Testing  effectiveness of control law (25) to suppress  spiral waves  with small
number  of  sensors  (2 equidistant in x-direction ( $y = 0.735L_y$ ) and   a grid of 25x6
actuators ( 25 equidistant   situated at  points 1,5,..,100 in x-direction and 6  at points

110,140,…,260  in y-direction, see Fig.6).  Control  start at t=10.   Other parameters are as in Figure 1.  Figure shows  snapshots at times  t=400, t=800,….

**Figure 5.** Testing  effectiveness of control law (25) to suppress  spiral waves  with  small number  of  sensors  (2 equidistant  sensors in x-direction ( $y = 0.735L_y$ ) and a grid of 50x6 actuators (50 equidistant  actuators  situated at  points 1,3,..,100 in x-direction and 6 at points 110,140,…, 260  in y-direction).  Control  start at t=10.   Other parameters are as in Figure 1.  Figure shows  snapshots at times  t=400, t=800,….

**Figure 6.**  Location of sensors ( ● ) and actuators (*) for control used in Figure 4.

**Figure 7.  Heterogeneous media:** Spiral waves in the (y,x) heterogeneous plane of the uncontrolled  dimensionless  transmembrane  potential  $V(x, y, t)$  (Eqns.1-3, $N_x = 200$ points,   $N_y = 400$ points,   parameters as in Fig.1 except   $\alpha$ ; $\alpha = 0.21$  in  the domain  $0.01L_x \leq x \leq 0.5L_x$, $0.275L_y \leq y \leq 0.775L_y$ during time $1 \leq t \leq 200$ as opposed to $\alpha = 0.1$ in Fig. 1. Figures shows  snapshots at times  t=100, t=200,….

**Figure 8.  Control in heterogeneous media:** Testing  effectiveness of control law (25) to suppress  spiral waves for conditions in Fig. 7, with  5 equidistant  sensors in x-direction ( $y = 0.735L_y$ ) and a grid of 50x6 actuators (50 equidistant  actuators  situated at  points 1,3,..,100 in x-direction and 6 at points 110,140,…, 260  in y-direction).  Control  start at t=10.   Other parameters are as in Figure 1.  Figure  shows  snapshots at times  t=200, t=400,….

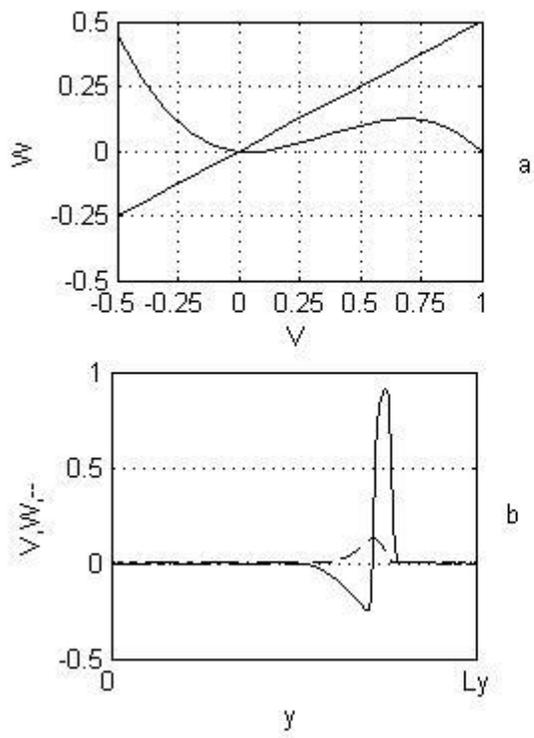

a

b

Figure 1

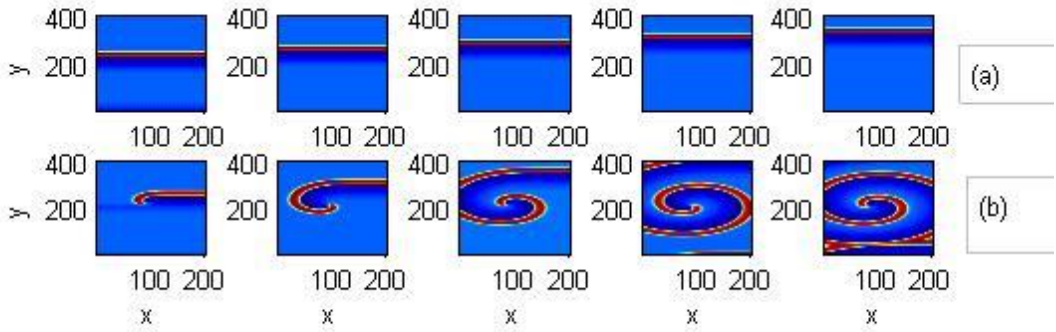

Figure 2

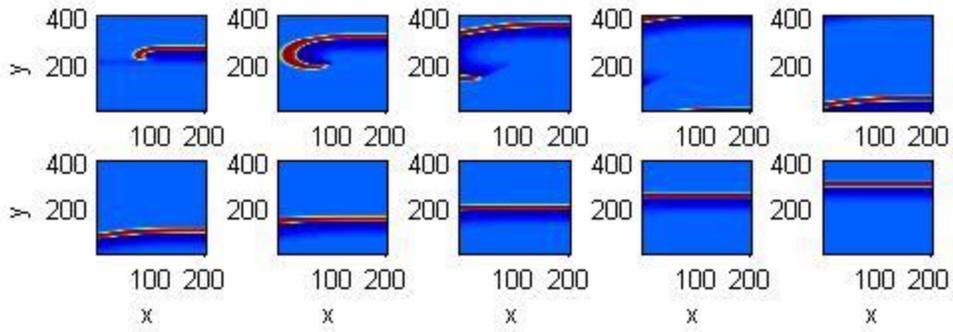

Figure 3

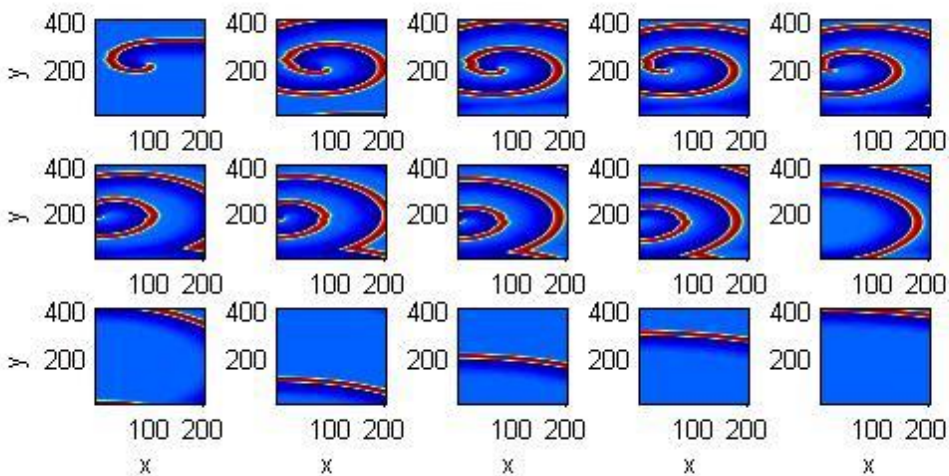

Figure 4

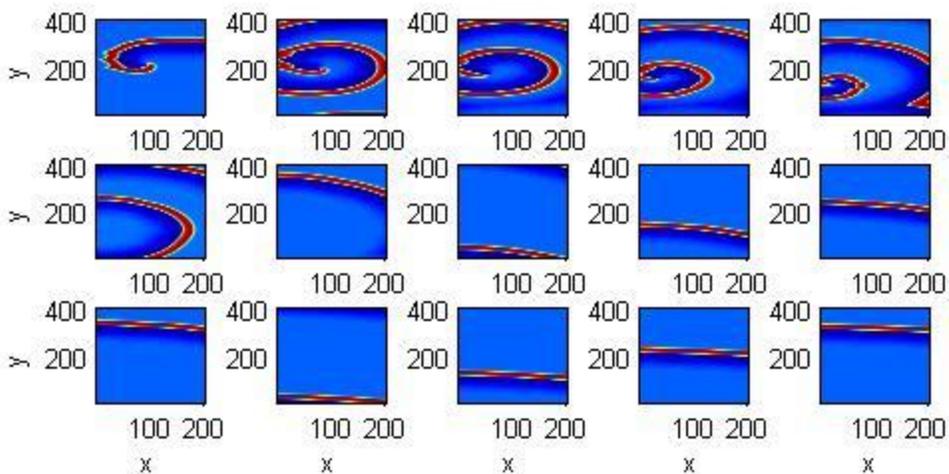

Figure 5

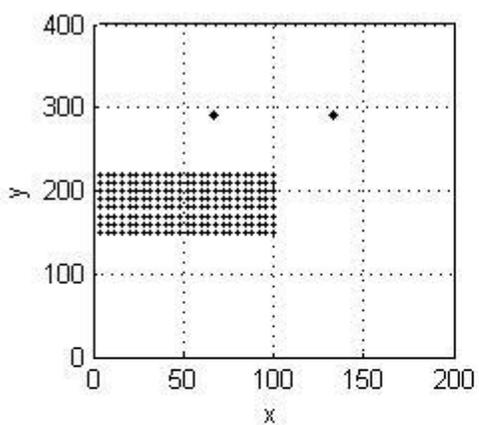

Figure 6

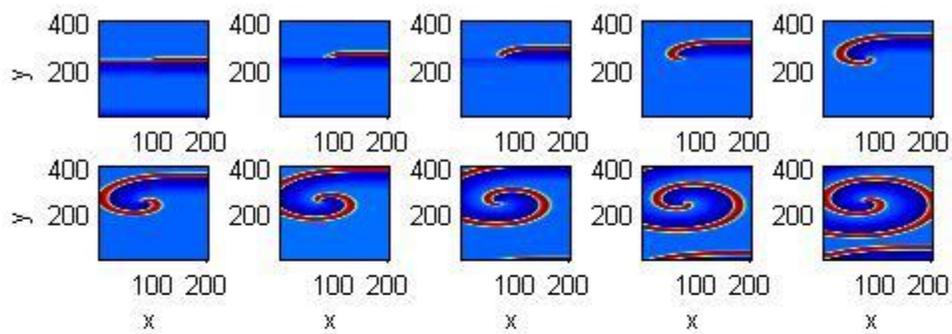

Figure 7

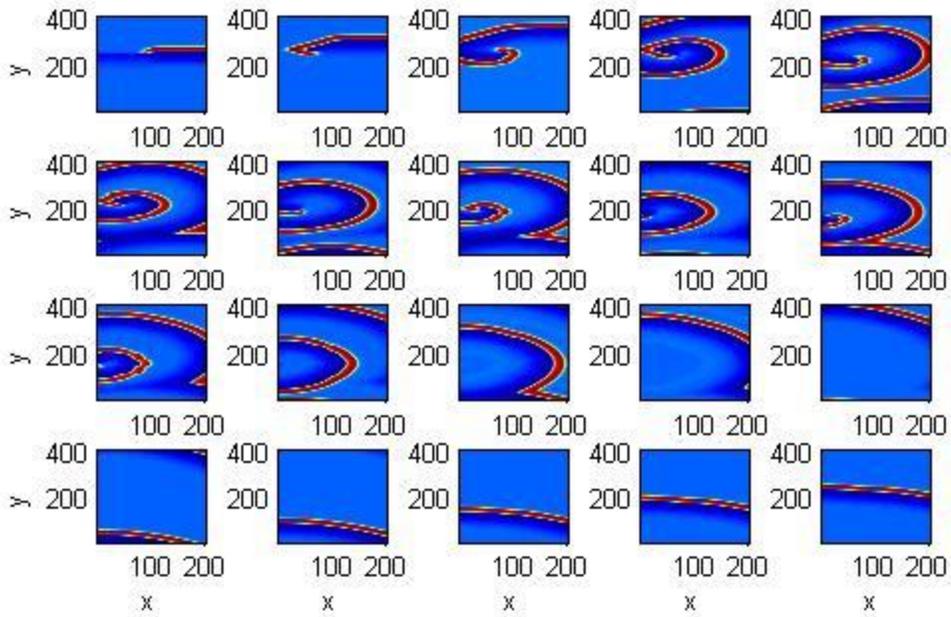

Figure 8